\documentclass[11pt]{article}
\usepackage{amssymb}
\makeatletter\@addtoreset {equation}{section}\makeatother

\newtheorem{theorem}{Theorem}
\newtheorem{lemma}{Lemma}
\newtheorem{remark}{Remark}
\newtheorem{example}{Example}
\newtheorem{corollary}{Corollary}
\newtheorem{definition}{Definition}
\newtheorem{assumption}{Assumption}
\newtheorem{proposition}{Proposition}

\newenvironment{proof}{
    \noindent {\it Proof.}}{\hfill$\Box$
}
\newenvironment{proof1}{
    \noindent {\it Proof }}{\hfill$\Box$
}

\usepackage[dvips]{epsfig}
\usepackage{graphicx}

\begin{document}

\title{\bf Justification of the coupled-mode approximation \\
for a nonlinear elliptic problem with a periodic potential}

\author{Dmitry Pelinovsky\footnote{\small On leave from Department of Mathematics, McMaster
University, Hamilton, Ontario, Canada, L8S 4K1} and Guido Schneider \\
{\small Institut f\"{u}r Analysis, Dynamik und
Modellierung Fakult\"{a}t f\"{u}r Mathematik und Physik,} \\
{\small Universit\"{a}t Stuttgart, Pfaffenwaldring 57, D-70569
Stuttgart, Germany} }

\date{\today}
\maketitle

\begin{abstract}
Coupled-mode systems are used in physical literature to simplify
the nonlinear Maxwell and Gross-Pitaevskii equations with a small
periodic potential and to approximate localized solutions called
gap solitons by analytical expressions involving hyperbolic
functions. We justify the use of the one-dimensional stationary
coupled-mode system for a relevant elliptic problem by employing
the method of Lyapunov--Schmidt reductions in Fourier space.
In particular, existence of periodic/anti-periodic and decaying
solutions is proved and the error terms are controlled in suitable
norms. The use of multi-dimensional stationary coupled-mode
systems is justified for analysis of bifurcations of
periodic/anti-periodic solutions in a small multi-dimensional
periodic potential.
\end{abstract}

\section{Introduction}

Gap solitons are localized stationary solutions of nonlinear
elliptic problems existing in the spectral gaps of the
Schr\"{o}dinger operator associated with a periodic potential. In
particular, gap solitons have been considered in two problems of
modern mathematical physics, the complex-valued Maxwell equation
\begin{equation}
\label{Maxwell} \nabla^2 E - \frac{n^2(x,|E|^2)}{c^2} E_{tt} = 0
\end{equation}
and the Gross-Pitaevskii equation
\begin{equation}
\label{GP} i E_t = - \nabla^2 E + V(x) E + \sigma |E|^2 E,
\end{equation}
where $E(x,t) : \mathbb{R}^N \times \mathbb{R} \mapsto \mathbb{C}$
and $\nabla^2 = \partial_{x_1}^2 + ... + \partial_{x_N}^2$. For
applications of the complex-valued Maxwell equation
(\ref{Maxwell}), $E$ is a complex amplitude of the electric field
vector, $c$ is the speed of light, and $n(x,|E|^2)$ is the
refractive index. The scalar equation (\ref{Maxwell}) is valid in
the space of one and two dimensions ($N = 1,2$) for so-called TE
modes but not in the space of three dimensions ($N = 3$), where
the system of Maxwell equations for a vector-valued function $E$
must be used \cite{Kuchment}. For applications of the
Gross--Pitaevskii equation (\ref{GP}), $E$ is the mean-field
amplitude, $\sigma$ is the scattering length, and $V(x)$ is the
trapping potential. The scalar equation (\ref{GP}) is the
mean-field model valid in the space of three dimensions ($N = 3$)
and it can be used in the space of one and two dimensions ($N =
1,2$) under additional assumptions \cite{bec_book}.

Stationary solutions of the Maxwell and GP equations are found
from the elliptic problem
\begin{equation}
\label{stationary} \nabla^2 U + \omega^2 U + \epsilon W(x) U =
\sigma |U|^2 U,
\end{equation}
where $U(x) : \mathbb{R}^N \mapsto \mathbb{C}$ and
($\omega^2,\epsilon,\sigma$) are parameters. The elliptic problem
(\ref{stationary}) is related to the GP equation (\ref{GP}) by an
exact reduction $E = U(x) e^{-i \omega^2 t}$, where the potential
$V(x)$ is represented by $V = -\epsilon W(x)$. The same problem is
related to the Maxwell equation (\ref{Maxwell}) with the
refractive index $n(x,|E|^2) = n_0^2 (1 + \mu W(x) + \nu |E|^2)$
by the exact reduction $E = U(x) e^{- i c \omega t/n_0}$, where
parameters are related by $\epsilon = \omega^2 \mu$ and $\sigma =
-\omega^2 \nu$.

Let us consider the elliptic problem (\ref{stationary}) with a
real-valued bounded potential $W(x)$, which is periodic in each
variable $x_j$, $\forall j$. The associated Schr\"{o}dinger
operator $L = - \nabla^2 - \epsilon W(x)$ is defined on
$C_0^{\infty}(\mathbb{R}^N)$ and is extended to a self-adjoint
operator which maps continuously $H^2(\mathbb{R}^N)$ to
$L^2(\mathbb{R}^N)$. Therefore, the spectrum $\sigma(L)$ is real.
Suppose that the absolutely continuous part of the spectrum of $L$
has a gap of finite size on the real spectral axis. For parameters
inside the spectral gap, localized solutions of the elliptic
problem (\ref{stationary}) were proved to exist in the relevant
variational problem \cite{Pankov}. (Earlier works on bifurcations
of gap solitons can be found in \cite{alama,heinz,stuart}.)
According to Theorem 1.1 of \cite{Pankov}, there exists a weak
solution $U(x)$ in $H^1(\mathbb{R}^N)$, which is (i) real-valued,
(ii) continuous on $x \in \mathbb{R}^N$ and (iii) decays
exponentially as $|x| \to \infty$.

We investigate more precise information on properties of the gap
soliton $U(x)$ by working with the asymptotic limit of small
$\epsilon$. When $\epsilon = 0$, the purely continuous spectrum of
$L_0 = - \nabla^2$ is non-negative and no finite gaps of
$\sigma(L_0)$ exist. However, when $N = 1$ and $0 \neq \epsilon
\ll 1$, narrow gaps of $\sigma(L)$ diverge from a sequence of
resonant points on a real axis and gap solitons may bifurcate
inside these narrow gaps. The coupled-mode system has been used in
the physical literature since the 1980s to characterize this
symmetry-breaking bifurcation of the spectrum $\sigma(L)$ and to
approximate one-dimensional gap solitons  of the problem
(\ref{stationary}) with $N = 1$ \cite{SS}. Therefore, our work
deals mainly with the case $N = 1$. We consider the potential
function $W(x)$ according to the following assumption.

\begin{assumption}
Let $W(x)$ be a smooth $2\pi$-periodic function with zero mean and
symmetry $W(-x) = W(x)$ on $x \in \mathbb{R}$. The Fourier series
representation of $W(x)$ is
\begin{equation}
\label{potential-series} W(x) = \sum\limits_{m \in \mathbb{Z}}
w_{2m} e^{i m x}, \quad \mbox{such that} \quad \sum_{m \in
\mathbb{Z}} (1 + m^2)^s |w_{2m}|^2 < \infty, \;\; \forall s >
\frac{1}{2},
\end{equation}
where $w_0 = 0$ and $w_{2m} = w_{-2m} = \bar{w}_{2m}$, $\forall m
\in \mathbb{N}$. \label{assumption-potential}
\end{assumption}

It will be clear from Proposition
\ref{proposition-singular-kernel} that the sequence of resonant
points of $\sigma(L)$ for $L = -\partial_x^2 - \epsilon W(x)$ is
located at $\omega = \omega_n = \frac{n}{2}$, $n \in \mathbb{Z}$,
so that a small gap in the spectrum $\sigma(L)$ bifurcates
generally from each point $\omega = \omega_n$, $n \in \mathbb{N}$
and a semi-infinite gap exists near $\omega = \omega_0 = 0$. A
formal asymptotic solution of the elliptic problem
(\ref{stationary}) near a resonant point $\omega = \omega_n$, $n
\in \mathbb{N}$, is given by
\begin{equation}
\label{formal_expansion} U(x) = \sqrt{\epsilon} \left[ a(\epsilon
x) e^{\frac{i n x}{2}} + b(\epsilon x) e^{-\frac{i n x}{2}} + {\rm
O}(\epsilon) \right], \qquad \omega^2 = \frac{n^2}{4} + \epsilon
\Omega + {\rm O}(\epsilon^2),
\end{equation}
where the vector $(a,b) : \mathbb{R} \mapsto \mathbb{C}^2$
satisfies the coupled-mode system with parameter $\Omega \in
\mathbb{R}$
\begin{equation}
\label{cme} \left\{ \begin{array}{c} i n a' + \Omega a + w_{2n} b
= \sigma (|a|^2 + 2 |b|^2) a, \\ -i n b' + \Omega b + w_{2n} a =
\sigma (2 |a|^2 + |b|^2) b, \end{array} \right.
\end{equation}
and the derivatives are taken with respect to $y = \epsilon x$.

The coupled-mode system (\ref{cme}) can be used for analysis of
bifurcations of periodic and anti-periodic solutions near narrow
gaps in the spectrum $\sigma(L)$. When the solutions $(a,b)$ of
system (\ref{cme}) are $y$-independent and the representation
(\ref{formal_expansion}) is used, the solution $\phi(x)$ of the
elliptic system (\ref{stationary}) at the leading order is
periodic in $x$ for even $n \in \mathbb{N}$ and anti-periodic for
odd $n \in \mathbb{N}$. It follows from system (\ref{cme}) that a
general family of $y$-independent small solutions $(a,b)$ near the
point $(0,0)$ has the form $a = c e^{i \theta}$ and $b = \pm c
e^{i \theta}$, where $c \in \mathbb{R}$, $\theta \in \mathbb{R}$,
and the nonlinear dispersion relation holds in the form
\begin{equation}
\label{nonlinear-dispersion-relation} \Omega \pm w_{2n} - 3 \sigma
c^2 = 0.
\end{equation}
The dispersion relation (\ref{nonlinear-dispersion-relation}) with
$c = 0$ shows that the eigenvalues for periodic/anti-periodic
solutions in the linear spectrum of the operator $L =
-\partial^2_x + \epsilon W(x)$ diverge from the point $\omega^2 =
\frac{n^2}{4}$ and result in a narrow gap in the spectrum
$\sigma(L)$ which lies in the interval
\begin{equation}
\label{narrow-gap} \frac{n^2}{4} - |\epsilon w_{2n}| < \omega^2 <
\frac{n^2}{4} + \epsilon |\epsilon w_{2n}|.
\end{equation}
The nonlinear dispersion relation
(\ref{nonlinear-dispersion-relation}) with $c \neq 0$ shows that
the nonlinear periodic/anti-periodic solutions bifurcate to the
right of the boundaries $\Omega = \pm w_{2n}$ for $\sigma = +1$
and to the left for $\sigma = -1$. We justify the persistence of
the leading-order results (\ref{nonlinear-dispersion-relation})
and (\ref{narrow-gap}) by using the method of Lyapunov--Schmidt
reductions in the discrete weighted space $l^2_s(\mathbb{Z}')$
equipped with the norm
\begin{equation}
\label{discrete-Sobolev-norm} {\bf U} \in l^2_s(\mathbb{Z}') :
\;\; \| {\bf U} \|^2_{l^2_s(\mathbb{Z}')} = \sum_{m \in
\mathbb{Z}'} \left( 1 + \frac{m^2}{4} \right)^s |U_m|^2 < \infty.
\end{equation}
Here $\mathbb{Z}'$ is a set of either even or odd numbers and
${\bf U}$ is a set of Fourier coefficients $\{ U_m\}_{m \in
\mathbb{Z}'}$ in the Fourier series
\begin{equation}
\label{solution-series} U(x) = \sqrt{\epsilon} \sum_{m \in
\mathbb{Z}'} U_m e^{\frac{i}{2} m x}, \qquad U_m = \frac{1}{2\pi
\sqrt{\epsilon}} \int_0^{2 \pi} U(x) e^{-\frac{i}{2} m x} dx,
\end{equation}
where the factor $\sqrt{\epsilon}$ is introduced for convenience.
The function $U(x)$ is periodic if the set $\mathbb{Z}'$ is even
and it is anti-periodic if $\mathbb{Z}'$ is odd. By the Sobolev
inequality, if ${\bf U} \in l^2_s(\mathbb{Z}')$ with $s >
\frac{1}{2}$, then the series (\ref{solution-series}) converges
absolutely and uniformly in the space of bounded continuous
functions $C_b^0(\mathbb{R})$ according to the bound
\begin{equation}
\label{potential-norm-bound} \sum_{m \in \mathbb{Z}'} |U_{m}| \leq
\sum_{m \in \mathbb{Z}'} (1+m^2)^s |U_m|^2 + \sum_{m \in
\mathbb{Z}'} \frac{1}{4(1 + m^2)^s} < \infty,
\end{equation}
where we have used the inequality $|a| |b| \leq |a|^2 +
\frac{1}{4} |b|^2$. Our main result on bifurcations of
periodic/anti-periodic solutions is summarized below.

\begin{theorem}
\label{proposition-LS-reductions}  Let Assumption
\ref{assumption-potential} be satisfied. Fix $n \in \mathbb{N}$,
such that $w_{2n} \neq 0$. Let $\omega^2 = \frac{n^2}{4} +
\epsilon \Omega$, where $\Omega \in \mathbb{R}$. The nonlinear
elliptic problem (\ref{stationary}) with $N = 1$ has a non-trivial
solution $U(x)$ in the form (\ref{solution-series}) with ${\bf U}
\in l^2_s(\mathbb{Z}')$ for any $s > \frac{1}{2}$ and sufficiently
small $\epsilon$ if and only if there exists a non-trivial
solution for $(a,b) \in \mathbb{C}^2$ of the bifurcation equations
\begin{equation}
\label{cme-local} \left\{ \begin{array}{c} \Omega a +
w_{2n} b - \sigma (|a|^2 + 2 |b|^2) a = \epsilon A_{\epsilon}(a,b), \\
\Omega b + w_{-2n} a - \sigma (2 |a|^2 + |b|^2) b = \epsilon
B_{\epsilon}(a,b),
\end{array} \right.
\end{equation}
where $A_{\epsilon}(a,b)$ and $B_{\epsilon}(a,b)$ are analytic
functions of $\epsilon$ near $\epsilon = 0$ satisfying the bounds
\begin{equation}
\label{bound-periodic-solution} \forall |\epsilon| < \epsilon_0,
\; \forall |a| + |b| < \delta : \quad |A_{\epsilon}(a,b)| \leq C_A
(|a| + |b|), \qquad |B_{\epsilon}(a,b)| \leq C_B (|a| + |b|),
\end{equation}
for sufficiently small $\epsilon_0 > 0$, fixed $\delta > 0$, and
some constants $C_A,C_B > 0$ which are independent of $\epsilon$
and depend on $\delta$. Moreover, $A_{\epsilon}(a,b) =
\bar{B}_{\epsilon}(b,a)$, $\forall (a,b) \in \mathbb{C}^2$ and
\begin{equation}
\label{representation-local} \forall |\epsilon| < \epsilon_0 :
\quad \left\| U(x) - \sqrt{\epsilon} \left( a e^{\frac{inx}{2}} +
b e^{\frac{-i n x}{2}} \right) \right\|_{C_b^0(\mathbb{R})} \leq C
\epsilon^{3/2}.
\end{equation}
for some $\epsilon$-independent constant $C > 0$.
\end{theorem}

\begin{corollary}
\label{corollary-LS-reductions} The coupled system
(\ref{cme-local}) admits a symmetry reduction $a = \bar{b}$, where
the value of $a \in \mathbb{C}$ satisfies the scalar equation
\begin{equation}
\label{cme-scalar} \Omega a + w_{2n} \bar{a} - 3 \sigma |a|^2 a =
\epsilon A_{\epsilon}(a,\bar{a}).
\end{equation}
Under the reduction, the solution $U(x)$ is real-valued.
\end{corollary}

The results of Theorem \ref{proposition-LS-reductions} and
Corollary \ref{corollary-LS-reductions} justify the use of the
$y$-independent coupled-mode system (\ref{cme}) for bifurcations
of periodic/anti-periodic solutions of the nonlinear elliptic
problem (\ref{stationary}) with $N = 1$. In particular, the only
non-zero solutions of the scalar equation (\ref{cme-scalar}) occur
for either $a \in \mathbb{R}$ or $a \in i \mathbb{R}$, when the
scalar equation (\ref{cme-scalar}) is reduced to the extended
nonlinear dispersion relation
\begin{equation}
\Omega \pm w_{2n} - 3 \sigma |a|^2 = \epsilon A_{\pm}(|a|),
\end{equation}
where $A_{\pm}(|a|) = \frac{1}{a} A_{\epsilon}(a,\bar{a})$ is a
bounded, real-valued error term for sufficiently small $\epsilon$
and finite value of $|a| \in \mathbb{R}$. Note that the values of
$A_{\pm}$ are real-valued due to the gauge invariance of the
coupled-mode system (\ref{cme-local}) inherited from the gauge
invariance of the elliptic problem (\ref{stationary}).

The newly formed gap (\ref{narrow-gap}) of the continuous spectrum
of $L = -\partial_x^2 + \epsilon W(x)$ corresponds to the interval
$|\Omega | < |w_{2n}|$. For instance, let $\sigma = +1$ and
$w_{2n} > 0$, then the localized solution of the coupled-mode
system (\ref{cme}) can be written in the exact form
\cite{ChPel,SS}
\begin{equation}
a(y) = \bar{b}(y) = \frac{\sqrt{2}}{\sqrt{3}} \frac{\sqrt{w_{2n}^2
- \Omega^2}}{\sqrt{w_{2n} - \Omega} \cosh(\kappa y) + i
\sqrt{w_{2n} + \Omega} \sinh(\kappa y)},\label{gap_soliton}
\end{equation}
where $\kappa = \frac{1}{n} \sqrt{w_{2n}^2 - \Omega^2}$. The exact
localized solution can easily be found for $\sigma = -1$ and
$w_{2n} < 0$. The trivial parameters of translations of solutions
in $y$ and $\arg(a)$ are set to zero in the explicit solution
(\ref{gap_soliton}), such that the functions $a(y)$ and $b(y)$
satisfy the constraints $a(y) = \bar{a}(-y) = \bar{b}(y)$.

\begin{definition}
\label{definition-orbit} The gap soliton of the coupled-mode
system (\ref{cme}) is said to be a {\em reversible homoclinic
orbit} if it decays to zero at infinity and satisfies the
constraints $a(y) = \bar{a}(-y) = \bar{b}(y)$.
\end{definition}

We justify the persistence of the gap soliton (\ref{gap_soliton})
of the coupled-mode system (\ref{cme}) in the nonlinear elliptic
problem (\ref{stationary}) by working with the Fourier transform
of $U(x)$
\begin{equation}
\label{solution-integral} U(x) =
\frac{\sqrt{\epsilon}}{\sqrt{2\pi}} \int_{\mathbb{R}} \hat{U}(k)
e^{i k x} dk, \qquad \hat{U}(k) = \frac{1}{\sqrt{2\pi \epsilon}}
\int_{\mathbb{R}} U(x) e^{-ikx} dx,
\end{equation}
where again the factor $\sqrt{\epsilon}$ is introduced for
convenience. We develop the method of Lyapunov--Schmidt reductions
in the continuous weighted space $L^1_q(\mathbb{R})$ equipped with
the norm
\begin{equation}
\label{norm-Banach} \hat{U} \in L^1_q(\mathbb{R}) : \;\; \|
\hat{U} \|_{L^1_q(\mathbb{R})} = \int_{\mathbb{R}} \left( 1 + k^2
\right)^{q/2} |\hat{U}(k)| dk  < \infty.
\end{equation}
By the Riemann--Lebesgue Lemma, if $\hat{U} \in L^1_q(\mathbb{R})$
for some $q \geq 0$, then the $n$-th derivative of the function
$U(x)$ is bounded and continuous for $0 \leq n \leq [q]$ and it
decays to zero at infinity, i.e., $U \in C^n_b(\mathbb{R})$ and
$\lim\limits_{|x| \to \infty} U^{(n)}(x) = 0$. Indeed, for any
$G(x) = U^{(n)}(x)$ and $\hat{G} \in L^1(\mathbb{R})$, it follows
that
\begin{equation}
\| G(x) \|_{C_b^0(\mathbb{R})} \leq C \| \hat{G}(k)
\|_{L^1(\mathbb{R})},
\end{equation}
for some $C > 0$. In addition, the Schwartz space is dense in
$L^1(\mathbb{R})$ such that there is a sequence $\{ \hat{G}_j(k)
\}_{j \in \mathbb{N}}$ in the Schwartz space which converges to
$\hat{G}(k)$ in $L^1$-norm, and therefore, there exists a sequence
$\{ G_j(x) \}_{j \in \mathbb{N}}$ which converges to $G(x)$ in
$C_b^0(\mathbb{R})$-norm, such that $\lim\limits_{|x| \to \infty}
G(x) = 0$.

Related to the coupled-mode system for $(a,b)$ in variable $y$, we
shall also use the Fourier transform for $(\hat{a},\hat{b})$ in
variable $p$, where $y = \epsilon x$ and $p = \frac{k}{\epsilon}$.
We note that the norm $L^1(\mathbb{R})$ is invariant as follows:
\begin{equation}
\label{scaling-transformation-invariance} \hat{A}(k) =
\frac{1}{\epsilon} \hat{a}\left( \frac{k}{\epsilon} \right) :
\qquad \| \hat{A} \|_{L^1(\mathbb{R})} = \| \hat{a}
\|_{L^1(\mathbb{R})}.
\end{equation}
This invariance explains the choice of the space
$L^1_q(\mathbb{R})$ in our analysis (see also \cite{SU}). Our main
result on the existence of gap soliton solutions is summarized
below.

\begin{theorem}
\label{theorem-main} Let Assumption \ref{assumption-potential} be
satisfied. Fix $n \in \mathbb{N}$, such that $w_{2n} \neq 0$. Let
$\omega = \frac{n^2}{4} + \epsilon \Omega$, such that $|\Omega| <
|w_{2n}|$. Let $a(y) = \bar{b}(y)$ be a reversible homoclinic
orbit of the coupled-mode system (\ref{cme}) in Definition
\ref{definition-orbit}. Then, the nonlinear elliptic problem
(\ref{stationary}) with $N = 1$ has a non-trivial solution $U(x)$
in the form (\ref{solution-integral}) with $\hat{U} \in
L^1_q(\mathbb{R})$ for any $q \geq 0$ and sufficiently small
$\epsilon$ such that
\begin{equation}
\label{bound-gap-soliton} \forall |\epsilon| < \epsilon_0 : \quad
\| U(x) - \sqrt{\epsilon} \left[ a(\epsilon x) e^{\frac{i n x}{2}}
+ \bar{a}(\epsilon x) e^{-\frac{i n x}{2}} \right]
\|_{C_b^0(\mathbb{R})} \leq C \epsilon^{5/6},
\end{equation}
for some sufficient small $\epsilon_0$ and $\epsilon$-independent
constant $C > 0$. Moreover, the solution $U(x)$ is real-valued,
continuous on $x \in \mathbb{R}$, and $\lim\limits_{|x| \to
\infty} U(x) = 0$.
\end{theorem}

The results of Theorem \ref{theorem-main} give more precise
information about gap solitons of the elliptic problem
(\ref{stationary}) with $N = 1$ compared to the general result in
Theorem 1.1 of \cite{Pankov}, since the leading-order
approximation of $U(x)$ is given by the exponentially decaying
solutions (\ref{gap_soliton}) of the coupled-mode system
(\ref{cme}). On the other hand, we do not prove in Theorem
\ref{theorem-main} that $U(x)$ decays exponentially as $|x| \to
\infty$.

Rigorous justification of the approximation
(\ref{formal_expansion}) and the time-dependent extensions of the
coupled-mode system (\ref{cme}) were developed in \cite{GWH} for
the system of cubic Maxwell equations and in \cite{SU} for the
Klein--Fock equation with quadratic nonlinearity. A bound on the
error terms was found in the Sobolev space $H^1(\mathbb{R})$ in
\cite{GWH} and in the space of bounded continuous functions
$C_b^0(\mathbb{R})$ in \cite{SU}. The bound is valid on a finite
interval of the time evolution, which depends on $\epsilon$. The
error is not controlled on the entire time interval $t \in
\mathbb{R}$ and, in particular, the formalism cannot be used for a
proof of persistence of the leading-order approximation
(\ref{formal_expansion}) and (\ref{gap_soliton}) for the
stationary solutions of the nonlinear elliptic problem
(\ref{stationary}). The results of our Theorem \ref{theorem-main}
are more precise than Theorem 1 of \cite{GWH} and Theorem 2.1 of
\cite{SU} in this sense, since the error bound of the
leading-order approximation is controlled independently of $t \in
\mathbb{R}$. In a similar context, the justification of the
nonlinear Schr\"{o}dinger equation for the nonlinear Klein--Gordon
equation with spatially periodic coefficients is reported in
\cite{KSU}.

The method of Lyapunov--Schmidt reductions for periodic solutions
was used in \cite{CW,venakidis}. The work \cite{CW} deals with a
two-dimensional lattice equation for the nonlinear wave equation
when eigenvalues of the relevant linearized operator accumulate
near the origin. In this case, the Nash--Moser Theorem must be
used for the infinite-dimensional part of the Lyapunov--Schmidt
decomposition. In our case, the linearized operator for the
one-dimensional lattice equation has eigenvalues bounded away of
the origin and the Implicit Function Theorem can be applied
without additional complications. This application of the
technique is similar to the one in \cite{venakidis}, which deals
with the periodic wave solutions in the system of coupled discrete
lattice equations. Other applications of the method for periodic
wave solutions in equations of fluid dynamics can be found in
\cite{Chen,CN}.

Persistence of modulated pulse solutions was considered in
\cite{GS01,GS05} in the context of the nonlinear Klein--Gordon
equations. (Earlier results on the same topics can be found in
\cite{AK,En}.) Methods of spatial dynamics were applied to a
relevant PDE problem for modulated pulse solutions, the
linearization of which possessed infinitely many eigenvalues on
the imaginary axis. The local center-stable manifold was
constructed for the nonlinear Klein--Gordon equations after
normal-form transformations and the pulse solutions were proved to
be localized along a finite spatial scale, while small oscillatory
tails occur generally beyond this spatial scale. In contrast to
these works, we will not reformulate the ODE problem as an
extended PDE problem and avoid the construction of the local
center-stable manifold. This simplification is only possible if
the variables of the time-dependent problems (\ref{Maxwell}) and
(\ref{GP}) can be separated and modulated pulse solutions are
described by the reduction to the elliptic problem
(\ref{stationary}). We note that the basic equations of
electrodynamics, such as the real-valued Maxwell equation (of the
Klein--Gordon type), would not support the separation of variables
and the modulated pulse solutions do not generally exist in the
real-valued Maxwell equation \cite{AK}.

The article is structured as follows. The proof of Theorem
\ref{proposition-LS-reductions} is given in Section 2, where the
technique of Lyapunov--Schmidt reductions in $l^2_s(\mathbb{Z}')$
is developed for bifurcations of periodic/anti-periodic solutions.
The proof of Theorem \ref{theorem-main} is given in Section 3,
where the technique of Lyapunov--Schmidt reductions in
$L^1_q(\mathbb{R})$ is extended for persistence of decaying
solutions. Section 4 discusses applications of similar methods for
the justification of multi-dimensional multi-component
coupled-mode systems with $N \geq 2$.

\section{Lyapunov--Schmidt reductions for periodic/anti-periodic solutions}

Let the potential $W(x)$ and the solution $U(x)$ to the elliptic
problem (\ref{stationary}) with $N = 1$ be expanded in the Fourier
series (\ref{potential-series}) and (\ref{solution-series})
respectively. By using the Fourier series, we convert the elliptic
problem (\ref{stationary}) with $N = 1$ in the space of bounded,
continuous, periodic/anti-periodic solutions $U \in
C_b^0(\mathbb{R})$ to a system of nonlinear difference equations
in the discrete Sobolev weighted space ${\bf U} \in
l^2_s(\mathbb{Z}')$ for some $s > \frac{1}{2}$. The nonlinear
difference equations are written in the explicit form
\begin{equation}
\label{difference-system} \left( \omega^2 - \frac{m^2}{4} \right)
U_m + \epsilon \sum\limits_{m_1 \in \mathbb{Z}'} w_{m-m_1} U_{m_1}
= \epsilon \sigma \sum\limits_{m_1 \in
\mathbb{Z}'}\sum\limits_{m_2 \in \mathbb{Z}'} U_{m_1}
\bar{U}_{-m_2} U_{m - m_1 - m_2}, \qquad \forall m \in
\mathbb{Z}',
\end{equation}
which can be casted in the equivalent matrix-vector form
\begin{equation}
\label{difference-singular} \left( {\cal L} + \epsilon {\cal W}
\right) {\bf U} = \epsilon \sigma {\bf N}({\bf U},\bar{\bf U},{\bf
U}).
\end{equation}
Here ${\bf U}$ is an element of the infinite-dimensional vector
space $l^2_s(\mathbb{Z}')$ with the norm
(\ref{discrete-Sobolev-norm}), elements of matrix operators ${\cal
L}$ and ${\cal W}$ are given by
$$
{\cal L}_{m,k} = \left( \omega^2 - \frac{m^2}{4} \right)
\delta_{m,k}, \quad {\cal W}_{m,k} = w_{m-k}, \qquad \forall (m,k)
\in \mathbb{Z}' \times \mathbb{Z}',
$$
and ${\bf N}({\bf U},\bar{\bf U},{\bf U}) = {\bf U} \star {\cal R}
\bar{\bf U} \star {\bf U}$ consists of the convolution operator
with the elements $({\bf U} \star {\bf V} )_m =
\sum\limits_{k\in\mathbb{Z}} U_k V_{m-k}$ and the inversion
operator with the elements $({\cal R} {\bf U})_m = U_{-m}$.

We shall verify that the nonlinear vector field associated with
the difference equations (\ref{difference-system}) is closed in
space ${\bf U} \in l^2_s(\mathbb{Z}')$ with $s > \frac{1}{2}$
(Lemma \ref{lemma-closed-iterations}). Working in this space, we
shall apply the Implicit Function Theorem in two cases $\omega
\neq \mathbb{R} \backslash \{ \frac{n}{2} \}_{n \in \mathbb{Z}}$
and $\omega = \omega_n = \frac{n}{2}$ for some $n \in \mathbb{N}$.
The first case is non-resonant and the Implicit Function Theorem
guarantees existence of the unique zero solution ${\bf U} = {\bf
0}$ of system (\ref{difference-singular}) near ${\bf U} = {\bf 0}$
and $\epsilon = 0$ (Lemma \ref{proposition-singular-kernel}). The
second case corresponds to a bifurcation of non-zero periodic or
anti-periodic solutions of system (\ref{difference-singular}) and
it is analyzed by using the Lyapunov--Schmidt decomposition. To
prove Theorem \ref{proposition-LS-reductions}, we will prove that
there exists a unique smooth map from the components
$(U_n,U_{-n})$ to the other components $U_m$, $\forall m \in
\mathbb{Z}' \backslash \{n,-n\}$. Projections to the components
$(U_n,U_{-n})$ yield the coupled-mode equations (\ref{cme-local}),
while the bounds on the remainder terms
(\ref{bound-periodic-solution}) follow from the bounds on the
vector ${\bf U}$ in space $l^2_s(\mathbb{Z}')$. Representation
(\ref{representation-local}) and symmetry reductions of Corollary
\ref{corollary-LS-reductions} follow from the technique of
Lyapunov--Schmidt reductions.

\begin{lemma}
\label{lemma-closed-iterations} Let ${\bf W} \in
l^2_s(\mathbb{Z})$ for all $s > \frac{1}{2}$. The vector fields
${\bf W} \star $ and ${\bf N}$ map elements of
$l^2_s(\mathbb{Z}')$ with $s
> \frac{1}{2}$ to elements of $l^2_s(\mathbb{Z}')$.
\end{lemma}

\begin{proof}
The space $l^2_s(\mathbb{Z})$ with $s > \frac{1}{2}$ forms a
Banach algebra. Therefore, there exists a constant $0 < C(s) <
\infty$ such that
\begin{equation}
\label{norm-algebra} \forall {\bf U},{\bf V} \in l^2_s(\mathbb{Z})
: \quad \| {\bf U} \star {\bf V} \|_{l^2_s(\mathbb{Z})} \leq C(s)
\| {\bf U} \|_{l^2_s(\mathbb{Z})} \| {\bf V}
\|_{l^2_s(\mathbb{Z})}, \quad \forall s > \frac{1}{2}.
\end{equation}
This property was proven in \cite{CW} and it is similar to the one
in the continuous $H^s(\mathbb{R})$ spaces. The maps ${\bf W}
\star {\bf U}$ and ${\bf U} \star {\cal R} \bar{\bf U} \star {\bf
U}$ act on ${\bf U} \in l^2_s(\mathbb{Z}')$, where $\mathbb{Z}'$
is a set of either even or odd numbers. Both convolution operators
transform a vector on $\mathbb{Z}'$ to a vector on $\mathbb{Z}'$.
Therefore, $l^2_s(\mathbb{Z}')$ forms a linear subspace in the
vector space $l^2(\mathbb{Z})$ with the same algebra property
(\ref{norm-algebra}). As a result, both ${\bf W} \star {\bf U}$
and ${\bf N}({\bf U},\bar{\bf U},{\bf U})$ map elements of
$l^2_s(\mathbb{Z}')$ with $s > \frac{1}{2}$ to elements of
$l^2_s(\mathbb{Z}')$.
\end{proof}

\begin{lemma}
\label{proposition-singular-kernel} Let ${\bf W} \in
l^2_s(\mathbb{Z})$ with $s > \frac{1}{2}$ and $\omega \in
\mathbb{R} \backslash \{ \frac{n}{2} \}_{n \in \mathbb{Z}}$ for $n
\in \mathbb{Z}$. The nonlinear lattice system
(\ref{difference-singular}) has a unique trivial solution ${\bf U}
= {\bf 0}$ in a local neighborhood of ${\bf U} = {\bf 0}$ and
$\epsilon = 0$.
\end{lemma}

\begin{proof}
If $\omega \in \mathbb{R} \backslash \{ \frac{n}{2} \}_{n \in
\mathbb{Z}}$, the operator ${\cal L}$ in system
(\ref{difference-singular}) is invertible and
$$
\| {\cal L}^{-1} \|_{l^2_s \mapsto l^2_s} \leq
\frac{1}{\min\limits_{m \in \mathbb{Z}'} |\omega^2 -
\frac{m^2}{4}|} = \rho_0 < \infty, \qquad \forall s \geq 0.
$$
It follows from the estimate (\ref{potential-norm-bound}) that
$\sum\limits_{k \in \mathbb{Z}'} |w_{m-k}| < \infty$ for any $m
\in \mathbb{Z}'$ and $s > \frac{1}{2}$. Therefore, the matrix
operator ${\cal W}$ is a relatively compact perturbation to ${\cal
L}$. By the perturbation theory \cite{Kato}, there exists an
$\epsilon$-independent constant $0 < \rho < \infty$ such that
$$
\| \left( {\cal L} + \epsilon {\cal W} \right)^{-1} \|_{l^2_s
\mapsto l^2_s} \leq \rho, \quad \forall s \geq 0,
$$
for sufficiently small $\epsilon$. By Lemma
\ref{lemma-closed-iterations}, the vector field on the right-hand
side of system (\ref{difference-singular}) is closed in
$l^2_s(\mathbb{Z}')$ for $s > \frac{1}{2}$. Moreover, it is
analytic with respect to ${\bf U} \in l^2_s(\mathbb{Z}')$ and
$\epsilon \in \mathbb{R}$. The zero solution ${\bf U} = {\bf 0}$
satisfies the nonlinear lattice system (\ref{difference-singular})
for any $\epsilon \in \mathbb{R}$. The Fr\'{e}chet derivative of
system (\ref{difference-singular}) at ${\bf U} = {\bf 0}$ (which
is just the operator ${\cal L} + \epsilon W$) has a continuous
bounded inverse for sufficiently small $\epsilon$. By the Implicit
Function Theorem, the zero solution ${\bf U} = {\bf 0}$ is unique
in a local neighborhood of ${\bf U} = {\bf 0}$ and $\epsilon = 0$.
\end{proof}

\begin{proof1}{\em of Theorem \ref{proposition-LS-reductions}}.
If $\omega = \frac{n}{2}$ for some $n \in \mathbb{N}$, the
operator ${\cal L}$ is singular with a two-dimensional kernel
$$
{\rm Ker}({\cal L}) = {\rm Span}\left( {\bf e}_{n}, {\bf e}_{-n}
\right) \subset l^2_s(\mathbb{Z}'),
$$
where ${\bf e}_n$ is a unit vector in $l^2_s(\mathbb{Z}')$. The
straightforward decomposition of $l^2_s(\mathbb{Z}') = {\rm
Ker}({\cal L}) \oplus {\rm Ker}({\cal L})^{\perp}$ is nothing but
the representation
\begin{equation}
\label{decomposition-U} {\bf U} = a {\bf e}_n + b {\bf e}_{-n} +
{\bf g},
\end{equation}
where
\begin{equation}
\label{constrained-space} {\bf g} \in {\rm Ker}({\cal L})^{\perp}
= \{ {\bf g} \in l^2_s(\mathbb{Z}') : \;\; g_n = g_{-n} = 0 \}.
\end{equation}
Let ${\cal P}$ be the projection operator from
$l^2_s(\mathbb{Z}')$ to ${\rm Ker}({\cal L})^{\perp}$ at $\omega^2
= \frac{n^2}{4}$. It is obvious that ${\cal P} {\cal L} {\cal P}$
is a non-singular operator at $\omega^2 = \frac{n^2}{4}$. By using
the same argument as in Lemma \ref{proposition-singular-kernel},
we obtain that there exists an $\epsilon$-independent constant $0
< \rho < \infty$, such that
\begin{equation}
\label{inverse-operator} \| \left( {\cal P} \left( {\cal L} +
\epsilon {\cal W} \right) {\cal P} \right)^{-1} \|_{l^2_s \mapsto
l^2_s} \leq \rho, \quad \forall s \geq 0
\end{equation}
for sufficiently small $\epsilon$. The inhomogeneous problem for
${\bf g}$ is written in the explicit form
\begin{eqnarray}
\nonumber m \neq \pm n : & \phantom{t} & \left( \frac{n^2 -
m^2}{4} + \epsilon \Omega \right) g_m + \epsilon \sum_{k \in
\mathbb{Z}'\backslash \{n,-n\}} w_{m-k} g_k  \\ & \phantom{t} &
\phantom{text} - \epsilon \sigma \sum_{m_1 \in
\mathbb{Z}'}\sum_{m_2 \in \mathbb{Z}'} U_{m_1} \bar{U}_{-m_2}
U_{m-m_1-m_2} = - \epsilon \left( a w_{m-n} + b w_{m+n} \right),
\label{inhomogeneous-problem}
\end{eqnarray}
where $U_m = g_m + a \delta_{m,n} + b \delta_{m,-n}$. By Lemma
\ref{lemma-closed-iterations}, the vector field of system
(\ref{inhomogeneous-problem}) is closed in $l^2_s(\mathbb{Z}')$
for $s > \frac{1}{2}$ and any $(a,b) \in \mathbb{C}^2$. Moreover,
it is analytic with respect to ${\bf g} \in l^2_s(\mathbb{Z}')$
for all $(a,b) \in \mathbb{C}^2$ and $\epsilon \in \mathbb{R}$. By
the bound (\ref{inverse-operator}) and the Implicit Function
Theorem there exists a unique trivial solution ${\bf g} = {\bf 0}$
of system (\ref{inhomogeneous-problem}) for $(a,b) = (0,0)$ and
any $\epsilon \in \mathbb{R}$. It is also obvious that the zero
solution exists for any $(a,b) \in \mathbb{C}^2$ and $\epsilon =
0$. For all $(a,b) \neq 0$, the Fr\'{e}chet derivative of system
(\ref{inhomogeneous-problem}) at ${\bf g} = {\bf 0}$ is different
from the matrix operator ${\cal P}({\cal L} + \epsilon W){\cal P}$
by the additional terms
$$
-\epsilon \sigma \left[ (|a|^2 + |b|^2) \delta_{m,k} + a \bar{b}
\delta_{m,k+2n} + \bar{a} b \delta_{m,k-2n} \right], \qquad
\forall (m,k) \in \mathbb{Z}' \times \mathbb{Z}'.
$$
For sufficiently small $\epsilon$ and finite $(a,b) \in
\mathbb{C}^2$, these terms change slightly the bound $\rho$ in
(\ref{inverse-operator}), such that the Fr\'{e}chet derivative
operator of system (\ref{inhomogeneous-problem}) at ${\bf g} =
{\bf 0}$ has a continuous bounded inverse for $|\epsilon| <
\epsilon_0$, where $\epsilon_0$ is sufficiently small. By the
Implicit Function Theorem, there exists a unique map ${\bf
G}_{\epsilon} : \mathbb{C}^2 \mapsto {\rm Ker}({\cal L})^{\perp}
\subset l^2_s(\mathbb{Z}')$, which is analytic in $\epsilon$ with
the properties ${\bf G}_{\epsilon}(0,0) = {\bf 0}$ and ${\bf
G}_0(a,b) = {\bf 0}$. Therefore, the map ${\bf G}_{\epsilon}$
admits the Taylor series expansion in $\epsilon$. The first term
of the Taylor series is
$$
{\bf G}_{\epsilon}(a,b) = \epsilon \left[ a {\bf g}_a + b {\bf
g}_b + a^2 \bar{b} {\bf g}_c + \bar{a} b^2 {\bf g}_d \right] +
{\rm O}(\epsilon^2),
$$
where non-zero components of vectors ${\bf g}_{a,b,c,d}$ in the
constrained space (\ref{constrained-space}) are
$$
({\bf g}_a)_m = \frac{4 w_{m-n}}{m^2 - n^2}, \quad ({\bf g}_b)_m =
\frac{4 w_{m+n}}{m^2 - n^2}, \quad ({\bf g}_c)_m = \frac{4 \sigma
\delta_{m,3n}}{m^2 - n^2}, \quad ({\bf g}_d)_m = \frac{4 \sigma
\delta_{m,-3n}}{m^2 - n^2}.
$$
Let $|a| + |b| < \delta$ and $\delta$ is fixed independently of
$\epsilon$. Due to the analyticity of ${\bf G}_{\epsilon}$ in
$\epsilon$, there exists an $\epsilon$-independent constant $C >
0$ such that
\begin{equation}
\label{map-bound} \forall |\epsilon| < \epsilon_0 : \quad \| {\bf
G}_{\epsilon}(a,b) \|_{l^2_s(\mathbb{Z}')} \leq \epsilon C \left(
|a| + |b| \right).
\end{equation}
The projection equations to the two-dimensional kernel of ${\cal
L}$ is found from system (\ref{difference-system}) at $m = \pm n$
in the explicit form
\begin{equation}
\begin{array}{c} (\Omega + w_0) a +
w_{2n} b - \sigma \sum\limits_{m_1 \in
\mathbb{Z}'}\sum\limits_{m_2 \in \mathbb{Z}'} U_{m_1}
\bar{U}_{-m_2} U_{n-m_1-m_2} = - \sum\limits_{k \in
\mathbb{Z}'\backslash \{n,-n\}} w_{n-k} g_k, \\
(\Omega + w_0) b + w_{-2n} a - \sigma \sum\limits_{m_1 \in
\mathbb{Z}'} \sum\limits_{m_2 \in \mathbb{Z}'} U_{m_1}
\bar{U}_{-m_2} U_{-n-m_1-m_2} = - \sum\limits_{k \in
\mathbb{Z}'\backslash \{n,-n\}} w_{-n-k} g_k
\end{array} ,  \label{homogeneous-system}
\end{equation}
where $U_m = a \delta_{m,n} + b \delta_{m,-n} + g_m$ and the map
${\bf g} = {\bf G}_{\epsilon}(a,b)$ is constructed above. At
$\epsilon = 0$, we obtain that ${\bf g} = {\bf 0}$ and
\begin{eqnarray*}
\sum\limits_{m_1 \in \mathbb{Z}'}\sum\limits_{m_2 \in \mathbb{Z}'}
U_{m_1} \bar{U}_{-m_2} U_{n-m_1-m_2} & = & \left( |a|^2 + 2 |b|^2
\right) a, \\
\sum\limits_{m_1 \in \mathbb{Z}'}\sum\limits_{m_2 \in \mathbb{Z}'}
U_{m_1} \bar{U}_{-m_2} U_{-n-m_1-m_2} & = & \left( 2 |a|^2 + |b|^2
\right) b.
\end{eqnarray*}
This explicit computation recovers the left-hand side of system
(\ref{cme-local}). The right-hand side is estimated from the bound
(\ref{map-bound}) on the map ${\bf G}_{\epsilon}(a,b)$ in
$l^2_s(\mathbb{Z}')$ with $s
> \frac{1}{2}$ to yield the bound (\ref{bound-periodic-solution}).

\noindent To prove the last assertion of Theorem
\ref{proposition-LS-reductions}, we shall prove that the map ${\bf
G}_{\epsilon}(a,b)$ has the symmetry $({\bf G}_{\epsilon})_m(a,b)
= (\bar{\bf G}_{\epsilon})_{-m}(b,a)$. Indeed, since $W(x)$ is
real-valued, its Fourier coefficients satisfy the constraint
$w_{-2m} = \bar{w}_{2m}$, $\forall m \in \mathbb{Z}$. The systems
of equations (\ref{inhomogeneous-problem}) and
(\ref{homogeneous-system}) are symmetric with respect to the
transformation $(a,b,g_m) \mapsto (b,a,\bar{g}_{-m})$. By
uniqueness of solutions of system (\ref{inhomogeneous-problem}) in
a local neighborhood of $\epsilon = 0$, we obtain that $({\bf
G}_{\epsilon})_m(a,b) = (\bar{\bf G}_{\epsilon})_{-m}(b,a)$,
$\forall m \in \mathbb{Z}'$. Then, it follows directly from system
(\ref{homogeneous-system}) that $A_{\epsilon}(a,b) =
\bar{B}_{\epsilon}(b,a)$.
\end{proof1}

\begin{proof1}{\em of Corollary \ref{corollary-LS-reductions}}.
Due to the property $A_{\epsilon}(a,b) = \bar{B}_{\epsilon}(b,a)$,
system (\ref{cme-local}) has the symmetry reduction $b = \bar{a}$,
which results in the scalar equation (\ref{cme-scalar}). The
vector ${\bf U}$ is given by the decomposition
(\ref{decomposition-U}) with $b = \bar{a}$ and $({\bf
G}_{\epsilon})_m(a,\bar{a}) = (\bar{\bf
G}_{\epsilon})_{-m}(\bar{a},a)$. Therefore, the solution $U(x)$
recovered from the Fourier series (\ref{solution-series}) is
real-valued.
\end{proof1}

\section{Lyapunov--Schmidt reductions for gap solitons}

Let the solution $U(x)$ to the elliptic problem (\ref{stationary})
with $N = 1$ be represented by the Fourier transform
(\ref{solution-integral}), while the potential $W(x)$ is given by
the Fourier series (\ref{potential-series}). The elliptic problem
(\ref{stationary}) with $N = 1$ is converted to the integral
advance-delay equation for the Fourier transform $\hat{U}(k)$:
\begin{equation}
\label{difference-system-gap} \left( \omega^2 - k^2 \right)
\hat{U}(k) + \epsilon \sum\limits_{m \in \mathbb{Z}} w_{2m}
\hat{U}(k-m) = \epsilon \sigma \int_{\mathbb{R}} \int_{\mathbb{R}}
\hat{U}(k_1) \hat{\bar{U}}(k_2) \hat{U}(k-k_1+k_2) dk_1 dk_2
\qquad \forall k \in \mathbb{R}.
\end{equation}
Working in the Fourier space $\hat{U} \in L^1_q(\mathbb{R})$,
where the vector field of the integral advance-delay equation
(\ref{difference-system-gap}) is closed (Lemma
\ref{lemma-closed-iterations-integral}), we decompose the solution
$\hat{U}(k)$ into three parts
\begin{equation}
\label{partition} \hat{U}(k) = \hat{U}_+(k)
\chi_{\mathbb{R}_+'}(k) + \hat{U}_-(k) \chi_{\mathbb{R}_-'}(k) +
\hat{U}_0(k) \chi_{\mathbb{R}_0'}(k),
\end{equation}
where $\chi_{[a,b]}(k)$ is a function of compact support (it is
$1$ on $k \in [a,b]$ and $0$ on $k \in \mathbb{R} \backslash
[a,b]$) and the intervals $\mathbb{R}'_+$, $\mathbb{R}'_-$ and
$\mathbb{R}_0'$ are
\begin{equation}
\mathbb{R}'_+ = \left[  \omega_n - \epsilon^{2/3}, \omega_n +
\epsilon^{2/3} \right], \quad  \mathbb{R}'_- = \left[ -\omega_n -
\epsilon^{2/3}, -\omega_n + \epsilon^{2/3} \right], \quad
\mathbb{R}_0' = \mathbb{R} \backslash (\mathbb{R}_+' \cup
\mathbb{R}_-').
\end{equation}
Here $\omega_n = \frac{n}{2}$ is a bifurcation value of $\omega$
and the components $\hat{U}_{\pm}(k)$ represent the largest part
of the solution $\hat{U}(k)$ near the resonant values $k = \pm
\omega_n$, which is approximated by the solution of the
coupled-mode system (\ref{cme}) in coordinates $y = \epsilon x$ in
physical space and $p = \frac{k}{\epsilon}$ in Fourier space. The
intervals surrounding the resonant values $k = \pm \omega_n$ have
small length $2 \epsilon^{2/3}$, where both the constant $c = 2$
and the scaling factor $r = \frac{2}{3}$ are fixed for
convenience. In fact, we could generalize all proofs for any
constant $c > 0$ and any scaling factor $\frac{1}{2} < r < 1$.

In order to prove Theorem \ref{theorem-main}, we shall apply the
method of Lyapunov--Schmidt reductions in space
$L^1_q(\mathbb{R})$ with $q \geq 0$. First, we prove the existence
of a unique smooth map from $(\hat{U}_+(k),\hat{U}_-(k))$ to
$\hat{U}_0(k)$ (Lemma \ref{lemma-existence-map}). The solutions
for the components $(\hat{U}_+(k),\hat{U}_-(k))$ are then
approximated by the suitable (exponentially decaying) solutions
$(\hat{a}(p),\hat{b}(p))$ of the coupled-mode system (\ref{cme})
rewritten in Fourier space (Lemma
\ref{lemma-estimate-gap-soliton}). The approximation yields the
desired bound (\ref{bound-gap-soliton}) for $\hat{U} \in
L^1_q(\mathbb{R})$ with $q \geq 0$. The reduction to the
real-valued solutions $U(x)$ becomes obvious from the
decomposition of the Lyapunov--Schmidt reduction method.
Continuity and decay conditions on $U(x)$ follow by the
Riemann--Lebesgue Lemma.

\begin{lemma}
\label{lemma-closed-iterations-integral} Let ${\bf W} \in
l^2_{s+q}(\mathbb{Z})$ for all $s > \frac{1}{2}$ and $q \geq 0$.
The vector field of the integral equation
(\ref{difference-system-gap}) maps elements of $L^1_q(\mathbb{R})$
with $q \geq 0$ to elements of $L^1_q(\mathbb{R})$.
\end{lemma}

\begin{proof}
The convolution sum in the integral equation
(\ref{difference-system-gap}) is closed due to the bound
$$
\forall \hat{U} \in L^1_q(\mathbb{R}), \; \forall {\bf W} \in
l^2_{s+q}(\mathbb{Z}) : \quad \left\| \sum_{m \in \mathbb{Z}}
w_{2m} \hat{U}(k-m) \right\|_{L^1_q(\mathbb{R})} \leq \| \hat{U}
\|_{L^1_q(\mathbb{R})} \|{\bf W} \|_{l^1_q(\mathbb{Z})} \leq \|
\hat{U} \|_{L^1_q(\mathbb{R})} \|{\bf W}
\|_{l^2_{s+q}(\mathbb{Z})}
$$
for any $q \geq 0$ and $s > \frac{1}{2}$, where the inequality
(\ref{potential-norm-bound}) has been used. The convolution
integral is closed due to the bound
$$
\forall \hat{U},\hat{V} \in L^1_q(\mathbb{R}) : \quad \left\|
\int_{\mathbb{R}} \hat{U}(k_1) \hat{V}(k-k_1) dk_1
\right\|_{L^1_q(\mathbb{R})} \leq C \| \hat{U}
\|_{L_q^1(\mathbb{R})} \| \hat{V} \|_{L_q^1(\mathbb{R})}
$$
for any $q \geq 0$ and some $C > 0$, which occur in the inequality
$$
1 + (k_1+k_2)^2 \leq C (1 + k_1^2)(1+k_2^2)
$$
for all $(k_1,k_2) \in \mathbb{R}^2$.
\end{proof}

\begin{lemma}
Let Assumption \ref{assumption-potential} be satisfied and $\omega
= \frac{n^2}{4} + \epsilon \Omega$, where $n \in \mathbb{N}$ and
$\Omega \in \mathbb{R}$. There exists a unique map
$\hat{U}_{\epsilon} : L^1_q(\mathbb{R}_+') \times
L^1_q(\mathbb{R}_-') \mapsto L^1_q(\mathbb{R}_0')$ for all $q \geq
0$, such that $\hat{U}_0(k) =
\hat{U}_{\epsilon}(\hat{U}_+,\hat{U}_-)$ and
\begin{equation}
\label{decay-bound} \forall |\epsilon| < \epsilon_0 : \quad \|
\hat{U}_0(k) \|_{L^1_q(\mathbb{R}_0')} \leq \epsilon^{1/3} C
\left( \|\hat{U}_+ \|_{L^1_q(\mathbb{R}_+')} + \|\hat{U}_-
\|_{L^1_q(\mathbb{R}_-')} \right),
\end{equation}
where $\epsilon_0$ is sufficiently small and the constant $C > 0$
is independent of $\epsilon$ and depends on $\delta$ in the bound
$\|\hat{U}_+ \|_{L^1_q(\mathbb{R}_+')} + \|\hat{U}_-
\|_{L^1_q(\mathbb{R}_-')} < \delta$ for a fixed
$\epsilon$-independent $\delta
> 0$. \label{lemma-existence-map}
\end{lemma}

\begin{proof}
We project the integral advance-delay equation
(\ref{difference-system-gap}) onto the interval $k \in
\mathbb{R}_0'$:
\begin{eqnarray*} \left( \frac{n^2}{4} + \epsilon \Omega
- k^2 \right) \hat{U}_0(k) + \epsilon \sum\limits_{m \in
\mathbb{Z}} w_{2 m} \chi_{\mathbb{R}_0'}(k) \hat{U}(k-m) =
\epsilon \sigma \chi_{\mathbb{R}_0'}(k) \int_{\mathbb{R}}
\int_{\mathbb{R}} \hat{U}(k_1) \hat{\bar{U}}(k_2)
\hat{U}(k-k_1+k_2) dk_1 dk_2,
\end{eqnarray*}
where $\hat{U}(k)$ is decomposed by the representation
(\ref{partition}). Since
$$
\min_{k \in \mathbb{R}_0'} \left| \frac{n^2}{4} - k^2 \right| \geq
C_n \epsilon^{2/3},
$$
for some $C_n > 0$, the linearized integral equation at $\epsilon
= 0$ is invertible such that
\begin{equation}
\label{decay-bound-local} \left\| \left( \frac{n^2}{4} - k^2
\right)^{-1} \right\|_{L^1_q(\mathbb{R}_0') \mapsto
L^1_{q}(\mathbb{R}_0')} \leq \frac{1}{C_n \epsilon^{2/3}}.
\end{equation}
The linearized integral equation on $\hat{U}_0(k)$ for $\epsilon
\neq 0$ is given on $k \in \mathbb{R}_0'$ by
\begin{eqnarray*}
\left( \frac{n^2}{4} + \epsilon \Omega - k^2 \right) \hat{U}_0(k)
+ \epsilon \sum\limits_{m \in \mathbb{Z}} w_{2m}
\chi_{\mathbb{R}_0'}(k) \hat{U}_0(k-m) \\ - \epsilon \sigma
\chi_{\mathbb{R}_0'}(k) \int_{\mathbb{R}_+'} \int_{\mathbb{R}_-'}
\left[ \hat{U}_+(k_1) \hat{\bar{U}}_-(k_2) + \hat{\bar{U}}_+(k_1)
\hat{U}_-(k_2) \right] \hat{U}_0(k + k_1+k_2) dk_1 dk_2 \\ -
\epsilon \sigma \chi_{\mathbb{R}_0'}(k) \int_{\mathbb{R}_+'}
\int_{\mathbb{R}_+'} \hat{U}_+(k_1) \hat{\bar{U}}_+(k_2)
\hat{U}_0(k-k_1+k_2) dk_1 dk_2 \\ -\epsilon \sigma
\chi_{\mathbb{R}_0'}(k) \int_{\mathbb{R}_-'} \int_{\mathbb{R}_-'}
\hat{U}_-(k_1) \hat{\bar{U}}_-(k_2) \hat{U}_0(k-k_1+k_2) dk_1
dk_2.
\end{eqnarray*}
Recall that $\epsilon^{2/3} \gg \epsilon$ for sufficiently small
$\epsilon$. If ${\bf W} \in l^2_{s+q}(\mathbb{Z})$ and
$\hat{U}_{\pm} \in L^1_q(\mathbb{R}_{\pm}')$ for all $s
> \frac{1}{2}$ and $q \geq 0$, the convolution sums and integrals are closed by Lemma
\ref{lemma-closed-iterations-integral}. Fix $\|\hat{U}_+
\|_{L^1_q(\mathbb{R}_+')} + \|\hat{U}_- \|_{L^1_q(\mathbb{R}_-')}
< \delta$ for some $\epsilon$-independent $\delta$. Then, the
linearized operator is continuously invertible for all $|\epsilon|
< \epsilon_0$. The integral equation is analytic in $\epsilon$ and
admits a unique trivial solution $\hat{U}_0(k) = 0$ if
$\hat{U}_{\pm}(k) = 0$ on $k \in \mathbb{R}_{\pm}'$ or if
$\epsilon = 0$. By the Implicit Function Theorem, there exists a
unique map $\hat{U}_{\epsilon} : L^1_q(\mathbb{R}_+') \times
L^1_q(\mathbb{R}_-') \mapsto L^1_q(\mathbb{R}_0')$ for $|\epsilon|
< \epsilon_0$. The map is analytic in $\epsilon$ near $\epsilon =
0$ with the properties $\hat{U}_0(\hat{U}_+,\hat{U}_-) = 0$ and
$\hat{U}_{\epsilon}(0,0) = 0$. Due to the analyticity of the map
$\hat{U}_{\epsilon}(\hat{U}_+,\hat{U}_-)$ in $\epsilon$ and the
bound (\ref{decay-bound-local}) on the inverse operator, the
solution $\hat{U}_0(k) = \hat{U}_{\epsilon}(\hat{U}_+,\hat{U}_-)$
satisfies the desired bound (\ref{decay-bound}).
\end{proof}

\begin{lemma}
\label{lemma-estimate-gap-soliton} Let Assumption
\ref{assumption-potential} be satisfied. Fix $n \in \mathbb{N}$,
such that $w_{2n} \neq 0$. Let $\omega = \frac{n^2}{4} + \epsilon
\Omega$, such that $|\Omega| < |w_{2n}|$. Let $a(y) = \bar{b}(y)$
be a reversible homoclinic orbit of the coupled-mode system
(\ref{cme}) in Definition \ref{definition-orbit}. Then, there
exists a solution of the integral equation
(\ref{difference-system-gap}), such that $\hat{U}_0(k) =
U_{\epsilon}(\hat{U}_+,\hat{U}_-)$ is given by Lemma
\ref{lemma-existence-map} and
\begin{eqnarray}
\label{closed-bound} \forall |\epsilon| < \epsilon_0 : \; \left\|
\hat{U}_+(k) - \frac{1}{\epsilon} \hat{a}\left(
\frac{k-\omega_n}{\epsilon} \right)
\right\|_{L^1_q(\mathbb{R}_+')} \leq C_a \epsilon^{1/3}, \;
\left\| \hat{U}_-(k) - \frac{1}{\epsilon} \hat{b}\left(
\frac{k+\omega_n}{\epsilon} \right)
\right\|_{L^1_q(\mathbb{R}_-')} \leq C_b \epsilon^{1/3},
\end{eqnarray}
for sufficiently small $\epsilon_0 > 0$ and $\epsilon$-independent
constants $C_a,C_b > 0$.
\end{lemma}

\begin{proof}
Let us use the scaling invariance
(\ref{scaling-transformation-invariance}) and map the intervals
$\mathbb{R}_{\pm}'$ for $\hat{U}_{\pm}(k)$ to the normalized
interval $\mathbb{R}_0 =
\left[-\epsilon^{-1/3},\epsilon^{-1/3}\right]$ for
\begin{equation}
\label{scaling-transformation} \hat{a}(p) = \epsilon
\hat{U}_+\left( \frac{k-\omega_n}{\epsilon} \right), \qquad
\hat{b}(p) = \epsilon \hat{U}_-\left( \frac{k +
\omega_n}{\epsilon} \right).
\end{equation}
The new functions $\hat{a}(p)$ and $\hat{b}(p)$ have a compact
support on $p \in \mathbb{R}_0$, while the norms $\| \hat{a}
\|_{L^1_q(\mathbb{R}_0)}$ and $\| \hat{b}
\|_{L^1_q(\mathbb{R}_0)}$ are equivalent to the norms $\|
\hat{U}_+ \|_{L^1_q(\mathbb{R}_+')}$ and $\| \hat{U}_-
\|_{L^1_q(\mathbb{R}_-')}$. Using the bound (\ref{decay-bound}),
we project the integral equation (\ref{difference-system-gap}) to
the system of two integral equations on $p \in \mathbb{R}_0$:
\begin{eqnarray}
\nonumber \left( \Omega + w_0 - n p \right) \hat{a}(p) + w_{2n}
\hat{b}(p) - \sigma \int_{\mathbb{R}_0} \int_{\mathbb{R}_0} \left[
\hat{a}(p_1) \hat{\bar{a}}(p_2) + \hat{b}(p_1)
\hat{\bar{b}}(p_2) \right] \hat{a}(p-p_1+p_2) dp_1 dp_2 \\
\label{integral-eq-1} - \sigma \int_{\mathbb{R}_0}
\int_{\mathbb{R}_0} \hat{a}(p_1) \hat{\bar{b}}(p_2) \hat{b}(p -
p_1+p_2) dp_1 dp_2 = \epsilon p^2 \hat{a}(p) + \epsilon^{1/3}
\hat{A}_{\epsilon}(\hat{a},\hat{b},\hat{U}_{\epsilon}(\hat{a},\hat{b})), \\
\nonumber \left( \Omega + w_0 + n p \right) \hat{b}(p) + w_{-2n}
\hat{a}(p) - \sigma \int_{\mathbb{R}_0} \int_{\mathbb{R}_0} \left[
\hat{a}(p_1) \hat{\bar{a}}(p_2) + \hat{b}(p_1)
\hat{\bar{b}}(p_2) \right] \hat{b}(p-p_1+p_2) dp_1 dp_2 \\
\label{integral-eq-2} - \sigma \int_{\mathbb{R}_0}
\int_{\mathbb{R}_0} \hat{b}(p_1) \hat{\bar{a}}(p_2) \hat{a}(p -
p_1+p_2) dp_1 dp_2 = \epsilon p^2 \hat{b}(p) + \epsilon^{1/3}
\hat{B}_{\epsilon}(\hat{a},\hat{b},\hat{U}_{\epsilon}(\hat{a},\hat{b})),
\end{eqnarray}
where $A_{\epsilon}$ and $B_{\epsilon}$ are computed from
$(\hat{a},\hat{b})$ and the map
$\hat{U}_{\epsilon}(\hat{a},\hat{b})$ of Lemma
\ref{lemma-existence-map}. When the right-hand side of system
(\ref{integral-eq-1})--(\ref{integral-eq-2}) is truncated and the
integration is extended to $p \in \mathbb{R}$, system
(\ref{integral-eq-1})--(\ref{integral-eq-2}) becomes the
coupled-mode system (\ref{cme}) rewritten after the Fourier
transform in $y$. If $w_{2n} \neq 0$ and $|\Omega| < |w_{2n}|$,
the coupled-mode system (\ref{cme}) has a reversible homoclinic
orbit $a(y) = \bar{b}(y)$. The Fourier transform $\hat{a}(p)$
decays exponentially as $|p| \to \infty$, such that the integrals
of the system (\ref{integral-eq-1})--(\ref{integral-eq-2}) on $p
\in \mathbb{R} \backslash \mathbb{R}_0$ are exponentially small in
$\epsilon$. Therefore, they can be moved to the right-hand side of
the system. In addition, the remainder terms $A_{\epsilon}$ and
$B_{\epsilon}$ are analytic in $\epsilon$ and controlled by the
bound (\ref{decay-bound}) on $(\hat{a}(p),\hat{b}(p))$ in the
space $L^1_q(\mathbb{R}_0)$ for all $q \geq 0$. The linear terms
$p^2 \hat{a}(p)$ and $p^2 \hat{b}(p)$ are also controlled by the
bound
\begin{equation}
\epsilon \| p^2 \hat{a}(p) \|_{L^1_q(\mathbb{R}_0)} \leq
\epsilon^{1/3} \| \hat{a}(p) \|_{L^1_q(\mathbb{R}_0)}, \qquad
\epsilon \| p^2 \hat{b}(p) \|_{L^1_q(\mathbb{R}_0)} \leq
\epsilon^{1/3} \| \hat{b}(p) \|_{L^1_q(\mathbb{R}_0)}.
\end{equation}
Therefore, system (\ref{integral-eq-1})--(\ref{integral-eq-2}) is
a perturbed coupled-mode system (\ref{cme}) in Fourier space with
the truncation error of the order ${\rm O}(\epsilon^{1/3})$
measured in space $L^1_q(\mathbb{R}_0)$. (Note that the system
(\ref{integral-eq-1})--(\ref{integral-eq-2}) is not closed on
$L^1_q(\mathbb{R})$ as the terms $p^2\hat{a}(p)$ and
$p^2\hat{b}(p)$ represent a singular second-order perturbation of
the first-order coupled-mode system.)

\noindent To prove the persistence of decaying solutions of the
coupled-mode system (\ref{cme}), we note that the linearized
differential operator associated to the coupled-mode system
(\ref{cme}) in the physical space is given by a self-adjoint
system of $4$-by-$4$ component Dirac operators:
\begin{equation}
\label{linearized-operator} \left[ \begin{array}{cccc} i n
\partial_y +  W_0 & - \sigma a^2
& w_{2n} - 2 \sigma a \bar{b} & -2 \sigma ab \\
- \sigma \bar{a}^2 & -i n \partial_y + W_0 &
- 2 \sigma \bar{a} \bar{b} & \bar{w}_{2n} - 2 \sigma a \bar{b} \\
\bar{w}_{2n} - 2 \sigma \bar{a} b & -2 \sigma a b & -i n
\partial_y + W_0 & - \sigma
b^2 \\ -2 \sigma \bar{a} \bar{b} & w_{2n} - 2 \sigma a \bar{b} &
-\sigma \bar{b}^2 & i n
\partial_y + W_0
\end{array} \right],
\end{equation}
where $W_0 = \Omega - 2 \sigma (|a|^2 + |b|^2)$. By Theorem 4.1
and Corollary 4.2 in \cite{ChPel}, the linearized operator
(\ref{linearized-operator}) is block-diagonalized into two
uncoupled $2$-by-$2$ Dirac operators, each has a one-dimensional
kernel. The two-dimensional kernel of the linearized operator
(\ref{linearized-operator}) is related to the translational
symmetries in $y$ and $\arg(a)$ with the eigenvectors
$[a'(y),b'(y),\bar{a}'(y),\bar{b}'(y)]^T$ and $[i a(y),i b(y),-i
\bar{a}(y),-i \bar{b}(y)]^T$. The zero eigenvalue of the
linearized operator (\ref{linearized-operator}) is bounded away
from the continuous spectrum and other eigenvalues on the real
axis \cite{ChPel}. The extended coupled-mode system given by the
system (\ref{integral-eq-1})--(\ref{integral-eq-2}) after the
Fourier transform is only solvable if the right-hand-side lies in
the range of the linearized operator (\ref{linearized-operator}).

\noindent The nonlinear elliptic problem (\ref{stationary}) with
real-valued symmetric potential $W(x) = W(-x) = \bar{W}(x)$ has
two symmetries: the gauge invariance $U(x) \to e^{i \alpha} U(x)$
for all $\alpha \in \mathbb{R}$ and the reversibility $U(x) \to
U(-x)$. The new system obtained after the Fourier transform
(\ref{solution-integral}) and the decomposition (\ref{partition})
inherits both symmetries, such that the extended coupled-mode
system is formulated in a constrained subspace orthogonal to the
kernel of the linearized operator (\ref{linearized-operator}). As
a result, the linearized operator is continuously invertible in
the constrained subspace of space $L^1_q(\mathbb{R}_0) \times
L^1_q(\mathbb{R}_0)$ for all $q \geq 0$. Truncation of the
integral terms introduces a small error in the remainder terms but
does not change the symmetries of the extended coupled-mode system
and does not alter the invertibility of the linearized operator.
By the Implicit Function Theorem, there exists a unique solution
of system (\ref{integral-eq-1})--(\ref{integral-eq-2}) for
$\hat{a}(p)$ and $\hat{b}(p)$ on $p \in \mathbb{R}_0$, which is
close to the reversible homoclinic orbit of the coupled-mode
system (\ref{cme}) in $L^1_q$-norm.
\end{proof}

\begin{remark}
{\rm If $w_{2n} > 0$ and $\Omega = w_{2n}$, the exact solution
(\ref{gap_soliton}) describes an algebraically decaying reversible
homoclinic orbit of the coupled-mode system (\ref{cme}). Since the
continuous spectrum of the linearized operator
(\ref{linearized-operator}) touches the zero eigenvalue in this
case, persistence of algebraically decaying reversible homoclinic
orbits can not be proved in Lemma
\ref{lemma-estimate-gap-soliton}.}
\end{remark}

\begin{remark}
{\rm The symmetry condition on the potential $W(-x) = W(x)$ in
Assumption \ref{assumption-potential} is important for the proof
of persistence of homoclinic orbits in the nonlinear elliptic
problem (\ref{stationary}) since it ensures that the set of
homoclinic orbits of the nonlinear problem (\ref{stationary})
includes the symmetric (reversible) homoclinic orbits $U(-x) =
U(x)$. Lemma \ref{lemma-estimate-gap-soliton} can not be proved if
the right-hand-side of the extended coupled-mode system is not in
the range of the linearized operator (\ref{linearized-operator}),
which may occur when a homoclinic orbit of the coupled-mode system
(\ref{cme}) is positioned arbitrarily with respect to the general
potential function $W(x)$. Non-persistence of such homoclinic
orbits is usually beyond all orders in the asymptotic expansion in
powers of $\epsilon$ and it is typical that the gap solitons
persist at two particular points on the period of the potential
$W(x)$ (see \cite{Pel} for details). In our paper, we avoid the
beyond-all-orders problem by imposing a symmetry condition on
$W(x)$ which is sufficient for existence of a reversible
homoclinic orbit which is centered at the point $x = 0$. }
\end{remark}

\begin{proof1}{\em of Theorem \ref{theorem-main}}.
When the solution $U(x)$ is represented by the Fourier transform
$\hat{U}(k)$, both scaling transformations of Lemma
\ref{lemma-existence-map} and \ref{lemma-estimate-gap-soliton} are
incorporated into the solution, and the bounds (\ref{decay-bound})
and (\ref{closed-bound}) are used, we obtain the bound
\begin{equation}
\forall |\epsilon| < \epsilon_0 : \quad \left\| \hat{U}(k) -
\frac{1}{\epsilon} \hat{a}\left(\frac{k - n/2}{\epsilon}\right) -
\frac{1}{\epsilon} \hat{b}\left(\frac{k + n/2}{\epsilon}\right)
\right\|_{L^1_q(\mathbb{R})} \leq C \epsilon^{1/3},
\end{equation}
which implies the desired bound (\ref{bound-gap-soliton}) in
original physical space. It remains to prove that the solution
$U(x)$ is real-valued. This property follows from the symmetry of
the map $U_0(k) = \hat{U}_{\epsilon}(\hat{U}_+,\hat{U}_-)$
constructed in Lemma \ref{lemma-existence-map} with respect to the
interchange of $ \hat{U}_+(k-\omega_n)$ and
$\bar{\hat{U}}_-(k+\omega_n)$ and the complex conjugation. As a
result, system (\ref{integral-eq-1})--(\ref{integral-eq-2}) has
the symmetry reduction $\hat{a}(p) = \hat{\bar{b}}(p)$, which is
satisfied by the solution of the truncated system. When the
partition (\ref{partition}) is substituted into the Fourier
transform (\ref{solution-integral}) with the symmetry
$\hat{U}_+(k-\omega_n) = \bar{\hat{U}}_-(k+\omega_n)$, the
resulting solution $U(x)$ is proved to be real-valued.
\end{proof1}

\section{Lyapunov--Schmidt reductions in multi-dimensional potentials}

Let us consider the elliptic problem (\ref{stationary}) in the
space of two dimensions ($N = 2$). Let the potential $W(x)$ be
periodic in both variables with the same normalized period, such
that
\begin{equation}
W(x_1 + 2\pi,x_2) = W(x_1,x_2+2\pi) = W(x_1,x_2), \qquad \forall
(x_1,x_2) \in \mathbb{R}^2.
\end{equation}
We shall justify the use of multi-component coupled-mode systems
for the analysis of bifurcations of two-dimensional
periodic/anti-periodic solutions of the elliptic problem
(\ref{stationary}). We use the Fourier series for the potential
$W(x)$ and the solution $U(x)$:
\begin{equation}
\label{potential-series-2D} W(x) = \sum\limits_{m \in
\mathbb{Z}^2} w_{2m} e^{i m \cdot x}, \qquad U(x) =
\sqrt{\epsilon} \sum_{m \in \mathbb{Z}_1' \times \mathbb{Z}'_2}
U_m e^{\frac{i}{2} m \cdot x},
\end{equation}
where $m \cdot x = m_1 x_1 + m_2 x_2$ and the sets $\mathbb{Z}'_1$
and $\mathbb{Z}'_2$ are even or odd if the solution $U(x)$ is
periodic or anti-periodic in the corresponding variable $x_1$ and
$x_2$. The elliptic problem (\ref{stationary}) with $N = 2$
transforms to a system of nonlinear difference equations, which is
similar to system (\ref{difference-system}):
\begin{equation}
\label{difference-system-2D} \left( \omega^2 - \frac{|m|^2}{4}
\right) U_m + \epsilon \sum\limits_{m_1 \in \mathbb{Z}'_1 \times
\mathbb{Z}'_2} w_{m-m_1} U_{m_1} = \epsilon \sigma
\sum\limits_{m_1 \in \mathbb{Z}'_1 \times \mathbb{Z}'_2}
\sum\limits_{m_2 \in \mathbb{Z}'_1 \times \mathbb{Z}'_2} U_{m_1}
\bar{U}_{-m_2} U_{m - m_1 - m_2},
\end{equation}
for all $m \in \mathbb{Z}'_1 \times \mathbb{Z}'_2$. The nonlinear
lattice system (\ref{difference-system-2D}) is closed in the space
$l^2_s(\mathbb{Z}'_1 \times \mathbb{Z}'_2)$ with $s > 1$ thanks to
the Banach algebra property:
\begin{equation}
\label{norm-algebra-2D} \forall {\bf U},{\bf V} \in
l^2_s(\mathbb{Z}^2) : \quad \| {\bf U} \star {\bf V}
\|_{l^2_s(\mathbb{Z}^2)} \leq C(s) \| {\bf U}
\|_{l^2_s(\mathbb{Z}^2)} \| {\bf V} \|_{l^2_s(\mathbb{Z}^2)},
\quad \forall s > 1,
\end{equation}
for some $C(s) > 0$. Under the same constraint $s > 1$, the double
Fourier series (\ref{potential-series-2D}) converges absolutely
and uniformly in $C^0_b(\mathbb{R}^2)$.

The system (\ref{difference-system-2D}) takes the same abstract
form (\ref{difference-singular}), where ${\bf U}$ is an element of
the vector space $l^2_s(\mathbb{Z}'_1 \times \mathbb{Z}'_2)$ with
$s > 1$. An extension of Lemma \ref{proposition-singular-kernel}
tells us that no non-trivial solution ${\bf U}$ of the system
(\ref{difference-system-2D}) exists in a local neighborhood of
${\bf U} = {\bf 0}$ and $\epsilon = 0$ unless $\omega = \omega_n =
\frac{|n|}{2}$, where $n = (n_1,n_2) \in \mathbb{Z}^2$ and $|n| =
\sqrt{n_1^2 + n_2^2}$. Bifurcations of non-trivial solutions occur
only in the resonant case $\omega = \omega_n$ and the number of
bifurcation equations (leading to the coupled-mode system) is
defined by the dimension of the resonant set $S_n$ in
\begin{equation}
\label{set-S-n} S_n = \left\{ m \in \mathbb{Z}'_1 \times
\mathbb{Z}'_2 : \quad |m|^2 = |n|^2 \right\}.
\end{equation}
Here again the set $\mathbb{Z}'_1$ is even/odd if $n_1$ is
even/odd and so is the set $\mathbb{Z}'_2$ with respect to $n_2$.

\begin{lemma}
\label{resonant-set} The set $S_n$ admits the following
properties:

\begin{itemize}
\item[(i)] $0 < {\rm Dim}(S_n) < \infty$.

\item[(ii)] If ${\bf n} = {\bf 0}$, the zero solution ${\bf m} =
{\bf 0}$ is unique.

\item[(iii)] If ${\bf n} = (n_1,0)$, then ${\rm Dim}(S_n) \geq 2$
if $n_1$ is odd and ${\rm Dim}(S_n) \geq 4$ if $n_1$ is even.

\item[(iv)] If ${\bf n} = (n_1,n_2) \in \mathbb{N}^2$, then ${\rm
Dim}(S_n) \geq 4$ if $n_1 - n_2$ is odd and ${\rm Dim}(S_n) \geq
8$ if $n_1-n_2$ is even and non-zero.
\end{itemize}
\end{lemma}

\begin{proof}
(i) follows from the bound $|m|^2 < \infty$ on the space of
integers and from the existence of the solution $m = n$. (ii) is
obvious from $|m|^2 = 0$. (iii) follows from the existence of
particular solutions $(\pm n_1,0)$ and $(0,\pm n_1)$ of $|m|^2 =
n_1^2$ (if $n_1$ is odd, the solutions $(0,\pm n_1)$ do not belong
to the space $\mathbb{Z}'_2$ of even numbers). (iv) follows from
the existence of particular solutions $(\pm n_1,\pm n_2)$ and
$(\pm n_2, \pm n_1)$ of $|m|^2 = n_1^2 + n_2^2$ (if $n_1-n_2$ is
odd, the solutions $(\pm n_2,\pm n_1)$ do not belong to the space
$\mathbb{Z}'_1 \times \mathbb{Z}'_2$ of the opposite parities and
if $n_1 = n_2$, the solutions $(\pm n_2, \pm n_1)$ are not
different from $(\pm n_1, \pm n_2)$).
\end{proof}

\begin{proposition}
\label{proposition-LS-reductions-2D} Let ${\bf W} \in
l^2_s(\mathbb{Z}^2)$ for all $s > 1$ and $\omega^2 =
\frac{|n|^2}{4} + \epsilon \Omega$ for some $n \in \mathbb{N}^2$
and $\Omega \in \mathbb{R}$. Let the set $S_n$ be defined by
(\ref{set-S-n}) with $d_S = {\rm Dim}(S_n)$. The nonlinear lattice
system (\ref{difference-system-2D}) has a non-trivial solution
${\bf U} \in l^2_s(\mathbb{Z}'_1 \times \mathbb{Z}'_2)$ for all $s
> 1$ and sufficiently small $\epsilon$ if and only if there exists
a non-trivial solution for ${\bf a} \in \mathbb{C}^{d_S}$ of the
bifurcation equations
\begin{equation}
\label{cme-local-2D} \Omega a_{j_m} + \sum_{m_1 \in S_n} w_{m -
m_1} a_{j_{m_1}} - \sigma \sum_{m_1 \in S_n} \sum_{-m_2 \in S_n'}
a_{j_{m_1}} \bar{a}_{j_{-m_2}} a_{j_{m - m_1 - m_2}} = \epsilon
A_{j_m,\epsilon}({\bf a}), \quad \forall m \in S_n,
\end{equation}
where $j_m$ is an index of $m$ in the set $S_n$, the set $S_n'$ is
a subset of $S_n$, such that  $m - m_1 - m_2 \in S_n$, ${\bf
A}_{\epsilon}({\bf a}) \in \mathbb{C}^{d_S}$ depends analytically
on $\epsilon$ near $\epsilon = 0$. Moreover, there exists constant
$C > 0$ which is independent of $\epsilon_0 > 0$ and depends on
$\delta > 0$ such that
\begin{equation}
\label{bound-periodic-solution-2D} \forall |\epsilon| <
\epsilon_0, \;\; \forall \|{\bf a} \|_{l^1(\mathbb{C}^{d_S})} <
\delta : \quad \|{\bf A}_{\epsilon}({\bf a})
\|_{l^1(\mathbb{C}^{d_S})} \leq C \|{\bf a}
\|_{l^1(\mathbb{C}^{d_S})},
\end{equation}
where $\epsilon_0$ is sufficiently small and $\delta$ is fixed
independently of $\epsilon_0$.
\end{proposition}

\begin{proof}
The proof repeats the proof of Theorem
\ref{proposition-LS-reductions} due to the fact that $|m|^2 -
|n|^2$ are bounded away from zero for $m \in \mathbb{Z}'_1
\times\mathbb{Z}'_2 \backslash S_n$. The Lyapunov--Schmidt
reductions are performed after the decomposition ${\bf U} =
\sum_{m \in S_n} a_{j_m} {\bf e}_m + {\bf g}$, where ${\bf g} \in
{\rm Ker}({\cal L})^{\perp}$.
\end{proof}

\begin{example}
{\rm The resonant value $\omega_{(0,0)} = 0$ corresponds to the
single-mode bifurcation, like in the one-dimensional problem
($N=1$). The next resonant value $\omega_{(1,0)} = \omega_{(0,1)}
= \frac{1}{2}$ corresponds to the two-mode bifurcation, which has
the same coupled-mode equations as in the one-dimensional problem
($N=1$) due to the separation of the periodic Fourier series in
one variable and the anti-periodic Fourier series in the other
variable. Finally, the next resonant value $\omega_{(1,1)} =
\frac{1}{\sqrt{2}}$ gives the first example of the non-trivial
four-component coupled-mode equations in the space of two
dimensions $(N = 2)$. The coupled-mode equations
(\ref{cme-local-2D}) can be rewritten explicitly for the
components $(a_1,a_2,a_3,a_4)$ which corresponds to the Fourier
modes for the resonant set $S_{(1,1)} = \left\{ (1,1); (-1,-1);
(1,-1); (-1,1)\right\}$ at the selected order:
\begin{eqnarray*}
(\Omega + w_{0,0}) a_1 + w_{2,2} a_2 + w_{0,2} a_3 + w_{2,0} a_4 &
= & \sigma \left( (|a_1|^2 + 2 |a_2|^2 + 2 |a_3|^2 + 2 |a_4|^2) a_1
+ 2 \bar{a}_2 a_3 a_4 \right), \\
(\Omega + w_{0,0}) a_2 + w_{-2,-2} a_1 + w_{-2,0} a_3 + w_{0,-2}
a_4 & = & \sigma \left( (2 |a_1|^2 + |a_2|^2 + 2 |a_3|^2 + 2
|a_4|^2) a_2 + 2 \bar{a}_1 a_3 a_4 \right), \\
(\Omega + w_{0,0}) a_3 + w_{2,-2} a_4 + w_{0,-2} a_1 + w_{2,0} a_2
& = & \sigma \left( (2 |a_1|^2 + 2 |a_2|^2 + |a_3|^2 + 2 |a_4|^2) a_3
+ 2 \bar{a}_4 a_1 a_2 \right), \\
(\Omega + w_{0,0}) a_4 + w_{-2,2} a_3 + w_{-2,0} a_1 + w_{0,2} a_2
& = & \sigma \left( (2 |a_1|^2 + 2 |a_2|^2 + 2 |a_3|^2 + |a_4|^2)
a_4 + 2 \bar{a}_3 a_1 a_2 \right).
\end{eqnarray*}
This system (with the derivative terms in $y_1 = \epsilon x_1$ and
$y_2 = \epsilon x_2$) was derived in \cite{AP05} by using
asymptotic multi-scale expansions. Higher-order resonances for
$\omega_n$ with larger values of $n \in \mathbb{N}^2$ may involve
more than four components in the coupled-mode equations
(\ref{cme-local-2D}), and the count of $d_{S}$ versus $n$ is not
available in general. }
\end{example}

\begin{remark}
{\rm Additional resonances were considered in \cite{AP05}, which
correspond to an oblique propagation of the resonant Fourier
modes, e.g. $e^{ \frac{i}{2} p x_1}$ and $e^{\frac{i}{2} ((p +
2m_1) x_1 + 2 m_2 x_2)}$ with $p = -(m_1^2 + m_2^2)/m_1 \notin
\mathbb{Z}$. These resonances can be incorporated in the present
analysis by using the transformation $U(x) = e^{\frac{i}{2} p x_1}
\tilde{U}(x)$. }
\end{remark}

\begin{remark}
\label{remark-non-existence} {\rm Existence of two-dimensional ($N
= 2$) gap soliton solutions can not be proved with the approach of
Section 3 for small values of $\epsilon$ when $\omega$ is close to
$\omega_n$. Indeed, if the solution $\hat{U}(k)$ is split into a
finite number of parts compactly supported near the points of
resonances $(k_1,k_2) = \frac{1}{2} (n_1,n_2)$ and the remainder
part $\hat{U}_0(\mathbb{R})$, then the operator $|k|^2 -
\omega_n^2$ with $\omega_n^2 = \frac{|n|^2}{4}$ is not invertible
in a neighborhood of the circle of the radius $|k| =
\frac{|n|}{2}$. Since only finitely many parts of the circle are
excluded from the compact support of $\hat{U}_0(k)$, the Implicit
Function Theorem can not be used to prove existence of the map
from the finitely many resonance parts of the solution
$\hat{U}(k)$ to the remainder part $\hat{U}_0(k)$. This obstacle
has a principal nature as it is related to a generic non-existence
of gap solitons in the systems without spectral gaps. Indeed, the
operator $L = -\nabla^2 - \epsilon W(x)$ has no gaps for
sufficiently small $\epsilon$ in the space of two dimensions ($N =
2$) \cite{Kuchment}.}
\end{remark}

\begin{remark}
{\rm The conclusions of Proposition
\ref{proposition-LS-reductions-2D} and Remark
\ref{remark-non-existence} can be extended to $N \geq 3$.}
\end{remark}

{\bf Acknowledgement.} D.P. thanks D. Agueev and W. Craig for
their help at the early stage of this project. The work of D.
Pelinovsky is supported by the Humboldt Research Foundation. The
work of G. Schneider is partially supported  by the
Graduiertenkolleg 1294 ``Analysis, simulation and design of
nano-technological processes'' sponsored by the Deutsche
Forschungsgemeinschaft (DFG) and the Land Baden-W\"{u}rttemberg.

\end{document}